 \documentclass{amsart}
 \usepackage[a4paper, margin=1.19in]{geometry}

\title{Morphisms (should be) everywhere}
\author{Attila Egri-Nagy$^1$ \and Mikl\'os Hoffmann$^2$}
\address{
$^1$Akita International University, Japan
  \and
  $^2$Eszterh\'azy K\'aroly Catholic University, Eger, Hungary.
}
\email{egri-nagy@aiu.ac.jp,\ hoffmann.miklos@uni-eszterhazy.hu}
\date{2025 June}
\usepackage{amsthm}
\theoremstyle{plain}
\usepackage{biblatex}
\usepackage[svgnames]{xcolor}
\usepackage[colorlinks,urlcolor=NavyBlue,citecolor=Maroon,linkcolor=DarkSlateGrey]{hyperref}
\usepackage{hyperref}

\usepackage{graphicx}%
\usepackage{multirow}%
\usepackage{amsmath,amssymb,amsfonts}%
\usepackage{amsthm}%
\usepackage{mathrsfs}%
\usepackage[title]{appendix}%
\usepackage{xcolor}%
\usepackage{textcomp}%
\usepackage{manyfoot}%
\usepackage{booktabs}%
\usepackage{todonotes}
\usepackage{tikz-cd}
\usepackage{tikz}
\usetikzlibrary{graphs}
\newtheorem{thesis}{Thesis}
\usepackage{chessboard}
\usepackage{array}

\newcommand{\morphism}[2]{\textbf{#1}$\rightarrow$\textbf{#2}}
\newcommand{\Morphism}[4]{(\textbf{#1}, \textbf{#2})$\rightarrow$(\textbf{#3}, \textbf{#4})}
\newcommand{\R}{\mathbb R}
\newcommand{\gap}{\vskip1ex}
\begin{document}


\begin{abstract}
  Morphisms, structure preserving maps, are everywhere in Mathematics as useful tools for thinking and problem solving, or as objects to study.
Here, we argue that the idea of operations being compatible across two domains goes beyond
its mathematical use: it is a fundamental mechanism of any intelligence.
We precisely define morphisms, distinguish between dynamic morphisms (on operations, binary relations) and static ones (on $n$-ary relations), and describe how a flexible and pluralistic use of morphisms can serve as a general framework for understanding and explanation in a wide variety of fields.
\end{abstract}


\maketitle

\tableofcontents

\section{Introduction}

\begin{quote}
 ``\ldots one should learn to think abstractly, because by doing so many
 philosophical difficulties simply disappear.'' \cite{gowers2002mathematics}
\end{quote}
What is understanding?
What is an explanation?
What are modeling and representation?
These are fascinating questions driving the philosophical discourse.
We do not claim to have the ultimate answers for these problems, but we think that a
disciplined application of a mechanism could make progress.
This mechanism likes to hide in plain sight, mostly behind the words `model', `representation' or `correspondence', or to disguise itself as mathematics-only tool.
Let us state upfront the main message.
\begin{thesis}
Finding compatible operations in two domains is not merely a tool in mathematics,
but a universal mechanism for any understanding, reasoning and explanation,
both for human and artificial intelligence.
\end{thesis}
These mechanisms are called \emph{morphisms} (i.e., structure-preserving maps),
and they could and should be used correctly and explicitly in all areas of knowledge
to make philosophical problems more precise and reduce the chance of confusion
arising from ambiguities.

In this paper we argue for the thesis in two different ways.

\begin{enumerate}
\item Morphisms are already in use in many fields beyond mathematics, but often
  not recognized as such.
\item When recognized, morphisms are frequently subject to misunderstanding and/or misuse, 
but with the correct application of this mechanism these misunderstandings can be reduced or eliminated.
\end{enumerate}
Establishing a morphism involves finding source and target domains and a correspondence between them
(this is often done implicitly in philosophical discourse), and paying special
attention to the compatible operations (this is what we emphasize).
Misunderstanding and misuse may arise from simple pairing versus compatible composable operations,
or from thinking about morphism in an unnecessarily rigorous way, e.g. only in
terms of totally-defined isomorphisms.

To avoid further misunderstandings, we precisely define morphisms and fix a
notation in Section
\ref{definition}.
We discuss the canonical examples in Section \ref{examples}.
In Section \ref{background} we go through some explicit and implicit uses of morphisms.
Section \ref{explanation} elaborates on how the idea of morphisms underlie
various aspects of fields related to intelligence, like explanation and understanding.
We conclude in Section \ref{conclusion} with a summary and a research plan.

\section{Morphisms: compatible operations}
\label{definition}

Here, we describe the general idea of morphisms with a degree of mathematical rigor sufficient for the current discourse.
This is a balancing act.
Being too specific could make morphisms look applicable only in Mathematics.
Being too general may also lead to misunderstandings.
Morphisms are not just some arbitrary similarities, there are clear defining rules.
For clarity, we mention upfront that
\begin{itemize}
  \item morphisms are dynamic (composition is at the core);
  \item they are not necessarily isomorphisms (they can be structure
forgetting);
\item they are not required to be functions (relations can have morphic
properties);
\item they come in pairs (we can always turn around a morphism, the reverse may cease to be a function, but it remains morphic);
\item and need not be totally defined (may not even require domains to
exist in full, only given by generators).
\end{itemize}
We will talk about relations on two levels: a morphism can be relational, and we can have a morphism of relations. 
To emphasize the dynamic nature of morphisms we introduce them in three steps,
in this particular order:
\begin{enumerate}
\item compatibility of composed operations;
\item operations as relations between things;
\item generating sets for morphisms.
\end{enumerate}

\subsection{Morphisms of Operations}

An \emph{operation} is a general term for doing something, some action.
A \emph{composition} of operations is just doing two operations in sequence, one after the
other.
The \emph{composite}, $a$ after $b$ denoted by $ab$,  is another operation.
We require that composition is \emph{associative}, or rather the operations we are interested in happen to be associative.
This means that $(ab)c=a(bc)$, thus we can group the operations as we like, the
result remains the same.
Consequently, the sequence $abc$ has a well-defined meaning.
Associativity states that there is nothing more than the linear order of $a,b$ and $c$, thus it captures an essential property of time \cite{wildbook}.

We call \emph{composable} the pairs of operations that can be connected in sequence, but not all of them can.
In computing, these constraints are enforced by a \emph{type system}.
Two computer programs are composable if the type of the output of the first
program is the same as the type of the input of the second one.
Two operations can be composable one way, but not necessarily the other way.
Making dough and then bake a bread form a composable pair of ``kitchen operations'', but trying to bake
flour, water, salt and yeast first and then trying to make a dough are not composable.

The importance of composition can be illustrated indirectly, to show what happens when it fails.
One can read a mathematical proof, being able to follow and verify every step to be correct, but still missing the understanding why the proof works.
The individual steps do not compose into a single composite.
A more direct demonstration is the immense modeling capability of category theory, where arrow composition is postulated in the definition of a category.

A \emph{domain of operations}  is defined by a set of operations serving as a
\emph{workspace} for thinking about a given field of knowledge.
This workspace is structured by composability relations, by its type system.
We may think, if it is helpful, that the operations are of the same kind, e.g.,
symbol manipulations or movements in physical space.
However, mathematically, we only need composability.
There is no need for ontological considerations.

We need two domains of operations.
Again, we may think that they are of different nature.
Indeed, the examples and applications will be like that.
Mathematically, this does not matter, they can be of the same kind.

A \emph{morphism} is a special relationship between two domains of
operations, from the \emph{source} to \emph{target}.
We will use the notation \morphism{Source}{Target} for naming morphisms by their domains.
For an operation in one domain we can find a corresponding operation in the
other one.
But not just any correspondence.
The operations have to be \emph{compatible with regard to composability}.
The idea of a morphism can be summarized by the rule:
\begin{equation}
  \label{eq:compatible_operations}
  \varphi(a\cdot b)=\varphi(a)\star\varphi(b),
\end{equation}
where $a,b$ are elements of the source domain.
Composition in the source domain is denoted by the $\cdot$ symbol, while in the target
domain by the $\star$ symbol.
Often in textbooks, the two compositions are not distinguished explicitly,
causing unnecessary confusion.
The narrative of the equation is easy to follow.
On the left side composition happens in the source domain, then we take the
composite operation to the target domain by $\varphi$.
On the right-hand side, we transfer the two operations to the target first, then
do the composition there.
The equation says that we get the same composite operation.
We require that this rule works for all composable pairs $a$ and $b$ in the
source domain.

We can express this with a commutative diagram (as in category theory).
The sets $S$ and $T$ denote the \emph{source} and the \emph{target} domains.
By $S\times S$ we mean a set of pairs of operations (direct
product), and $\varphi\times \varphi$ says that we apply $\varphi$ to both
elements of the ordered pair.
The diagram commutes
if different paths do the same thing.

\begin{center}
\begin{tikzcd}
  S\times S \ar[d,"\varphi\times \varphi"]
  \ar[r,"\cdot"]& S \ar[d,"\varphi"] \\
  T\times T \ar[r, "\star"]& T
\end{tikzcd}
\end{center}
Starting from $S\times S$ we can go right, meaning that we do the composition in
$S$ and then transfer the results to $T$.
On the level of individual elements, $a,b\in S$, this is $\varphi(a\cdot b)$.
We can also take the other route (down first, then right), transfer the
operations and do the composition in $T$.
This is $\varphi(a)\star \varphi(b)$.

What is $\varphi$? It can be a function $S\rightarrow T$.
It is important to emphasize, that it does not need to be a bijective function, leading to \emph{isomorphisms}.
More general functions can represent abstraction, lossy compression, yielding \emph{homomorphisms}.
Moreover, $\varphi$ can also be a
relation, which is a set-valued function $S\rightarrow 2^T$, where $2^T$ denotes
the power set of $T$.
For a \emph{relational morphism}, the morphic property is expressed by
$$\varphi(a)\varphi(b)\subseteq \varphi(ab).$$
Relations can be turned around easily, unlike functions, where swapping domains, in general, leads to a relation.
This feature was the original motivation in group theory  for defining relational morphisms (\cite{Wedderburn1941}),
and they are also one of the main tools for semigroups (\cite{QBook}).


\subsection{Morphisms of Binary Relations}

There is one operation of particular interest: tracing a binary relation.
When two things are related, we draw an arrow between them.
Moving along arrows is composable for head-to-tail sequences of arrows.
Through this operation we can define morphisms for binary relations.

We consider an individual relation as the most basic building block for our theory of morphisms.  
It can be represented
by a simple diagram.
\begin{center}
  \begin{tikzcd}
    x\ar[r,"a"] & y
  \end{tikzcd}
\end{center}
These relations are composable.
The composition $b\circ a$ reads as $b$ after $a$.
We can denote the movement from $x$ to $y$ by $xa$, so $xa=y$.
\begin{center}
  \begin{tikzcd}
    x\ar[r,"a"]\ar[rr, bend left,  "b\circ a"] & y \ar[r,"b"] & z
  \end{tikzcd}
\end{center}

There are numerous possible interpretations:  $x$ is $a$-related to $y$; $a$ takes
$x$ to $y$; $x$ moves to $y$ under the condition $a$.
We think that this elementary relation can be taken as a fundamental motif.
Traditionally, mathematics starts with sets, but `belonging to a set' is a special case of an elementary relation.

When talking about relations, we have two kinds of entities: the things and the
relations between them.
In graph theory we have \emph{nodes} and \emph{edges}, in category theory we
have \emph{objects} and \emph{arrows}.
Consequently, for defining a morphism we need to specify two functions (or relations).
We often use a geometrical language (rather metaphorically) for describing these two parts of the morphism: $\varphi_0$ maps the 0-dimensional points, the objects, and $\varphi_1$
maps the 1-dimensional arrows.
The condition of compatible operations is
$$\varphi_0(xa)=\varphi_0(x)\varphi_1(a).$$
The narrative of the formula is a little more involved.
We can stay in the source domain, and apply relation $a$ to object $x$, and then map the result to the target domain.
Or, we can map the object and the relation arrow to the target domain, and do the application there.
Either way, we should get the same object in the target domain.

Tracing a binary relation between objects and acting by an operation on objects are essentially the same.
Therefore, shifting the terminology from relation to operation is justified.
We will use the general notation for morphisms of this kind:
\gap
\begin{center}\Morphism{things}{operations}{other things}{other operations}.
\end{center}
\gap
Alternatively, using the geometric metaphor for more precision we can also write:
\gap
\begin{center}
  \Morphism{Source$_0$}{Source$_1$}{Target$_0$}{Target$_1$}.
\end{center}
\gap
This definition of morphisms is close to category theory.
We have graphs where all composable paths have dedicated arrows.
However, we do not require to have identity arrows on the objects, thus the domains we talk about are semigroupoids.

\subsection{Morphisms of $n$-ary relations -- no easy composition}

We can define morphisms for $n$-ary relations too.
A morphism between $S$ and $T$ is the function $\varphi: S\rightarrow T$, such that for all $n$-tuples $(x_1, \ldots, x_n)$ that are in relationship in $S$, the relation in $T$ holds for the $n$-tuple $(\varphi(x_1),\ldots,\varphi(x_2))$.

However, an $n$-ary relation in general does not have directionality like the
binary case does.
We can say that the $n$ elements of related, but relation does not specify any further roles, distinction, any further connections between them. 
If we want to impose direction on the connections, then we have two different possibilities already for $n=3$.
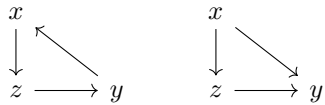
\begin{center}
  \begin{tikzcd}
    x \ar[d,dash] &  \\
    z \ar[r,dash]& y \ar[ul,dash]
  \end{tikzcd}
\hskip4em
  \begin{tikzcd}
    x \ar[d] &  \\
    z \ar[r]& y \ar[ul]
  \end{tikzcd}
\hskip2em
  \begin{tikzcd}
    x \ar[d] \ar[dr] &  \\
    z \ar[r]& y
  \end{tikzcd}

\end{center}
Consequently, we do not have composition naturally, not in a sequence.
Therefore, morphisms between $n$-ary relations can only transfer \emph{static}
structure (how things are related).
This limitation reduces the direct applicability of these morphisms.
We do not deny the usefulness of $n$-ary relations, but we claim that understanding is \emph{dynamic} even for static structures, i.e., it happens through a sequence of composable events.
When observing an unchanging object, we still have an order of looking at the parts.
This sequence defines movement, which can be matched to compatible movements in other domains.
In a way, thinking is movement.

\subsection{Generated Domains}

In mathematics, we consider the domains of morphisms to be complete, even for infinite ones.
In the wide range of applications we suggest, most morphisms are defined only for a subset of the domain.
However, some general inferential process often allows us to build the domain and the morphic relation on demand.
A subset of a domain from which we can potentially reconstruct the whole domain is called a \emph{generator set}.
For instance, zero, one, and the addition operator generate the set of natural numbers.
The $\frac{2\pi}{n}$ rotation generates all rotations of a regular $n$-gon.
The input symbols of a finite state automaton generate all possible computations it can perform.
A set of axioms and logical inference generates all true statements.

In practice, we need the morphism to be defined only for generator objects and operations in the source domain; the rest of it follows from the condition for compatible operations.
Often, the real question is how long the generation process compatibility holds.

\section{Canonical Examples}
\label{examples}

\subsection{The Cartographic Map -- The Basic idea}
In both scientific and everyday life, the most obvious form of appearance of morphism is a map - the map we follow to find something somewhere, the map from which we determine where we are in relation to something or somebody else.
Tracing lines on the map correspond to walking in town. In everyday practice it works well, as the two objects, the city and the drawing, are in very close correspondence.

In the source domain, for instance, the operations are to follow the road with your finger, choose the middle one at the fork, and checking that we have to turn left at the next fork.
In fact, we connect vertices on the road graph and compose these edges.
In reality, we do this with our car, driving from one corner to the other, composing the individual segments of roads. What we mean by the fact that the map is \emph{compatible with regard to composability} is the following: if we look at the map at the triple junction, then look up and choose the middle road with our car, then at the next corner we look at the map again, and then with our car we turn to the left, we get to the same place as if we had looked on the map to observe that first we have to take the middle road at the triple junction, then the left exit at the next junction, and finally we drive along this long road without even looking at the map.

There is also a use case for $n$-ary relations.
When being lost, we use the map to find the current location.
We look around and identify some landmark and note its relationships with other
features of the city.
Then we  can find some candidate location on the map and check whether the same
relationships hold.
This is matching static patterns, $n$-ary relations.
Note that if we keep moving while trying to find the location, to get different
perspectives, then we also use a morphism of operations.

If we accept the assumption that thinking in general is movement in an abstract space of thoughts (see Section \ref{sect:1000brains}), then using maps becomes more than being a prime example of morphisms.
Mapmaking can serve as a mechanism for explanations and understanding.

\subsection{The Cartographic Map -- Mathematical Background}

In the cartographic sense, the map, as a planar representation of (part of) the known world, is a very ancient invention. However, in a mathematical sense, it was only at the beginning of the 19th century that Gauss put the challenges and problems of mapmaking on a theoretical basis. There is no perfect map, and this caused just as much trouble (e.g. for sailors) in Gauss's time as before.

What do we mean by ``perfect map''? In the ideal case, the map, that is the \emph{source} domain represents the part to be depicted, the \emph{target} domain in an isometric manner, i.e. preserving all metric properties, such as distance, angle, area, to the extent of proportionality. Gauss proved that two surfaces can only be mapped onto each other isometrically if their curvatures are equal from point to point. That is, it is impossible to create an isometric planar map of the earth's surface, due to its varied curvature, hills and valleys. It is interesting that the otherwise very prolific Gauss considered this theorem as his main mathematical result, calling it Theorema Egregium.

Even if there is no ``perfect'' map, we still need a map. If the map is not accurate metrically, then we are faced with a choice: it could be absolutely accurate in some respects and very distorted in other respects, or it could be softly distorted from each metric point of view, but not very bad in any respect. There are many types of maps. An example of the former is an angle-preserving (or, from mathematical point of view, conformal) map, which, in turn, greatly spoils the distances. But most often we prefer the latter case, i.e. we use a so-called ``generally distorted'' map. In any case, however, we certainly expect a structure-preserving map from the \emph{source} domain to the \emph{target} domain, that is a morphism. If we see a three-way junction on the map, from which we have to choose the middle one, and then turn left to get to the store, we expect the same structure and operation to be present and to be successful in reality.

Another level of the mapping problem is when we want to gain knowledge about an entity that is theoretically or practically impossible to map as a whole, but some parts of it can still be mapped well, i.e. if not globally, there are local maps. In a mathematical sense, this area is the mathematics of manifolds introduced by Riemann, where the world of interest behaves locally in a Euclidean way (i.e. it can be mapped locally to $\R^n$), but not necessarily globally.

What do we mean by local mapping? There is a topological space $M$, and its subset $U$. A chart $(U, \varphi)$ is a (homeo)morphism from $U$ to an open subset of the Euclidean space $\R^n$. This chart allows us to locally understand $M$. For example, if a function $f : M \rightarrow \R$ is given, and we have a chart $(U, \varphi)$ on M, one could consider the composition $f \circ \varphi^{-1}$ which is a kind of explanation of the behaviour of the real-valued function on $M$ based on its behaviour on the chart.

But this situation needs some further considerations. Let $(V, \psi)$ be another chart on M, and suppose that U and V have some intersection. It would be essential to define the two maps in a way that those parts in common on $M$ have ``similar'' images in the two maps. The two maps should be consistent, that is, the two corresponding functions $f \circ \varphi^{-1}$ and $ f \circ \psi^{-1}$ should be linked in the sense that they can be reparametrized into one another as

  $${\displaystyle f\circ \varphi^{-1}={\big (}f\circ \psi^{-1}{\big )}\circ {\big (}\psi \circ \varphi^{-1}{\big )}.}$$

Summarizing our experience with this canonical example: even in this very ordinary modeling issue, map making, it can be seen that insisting on isomorphism is not tenable. 
Moreover, in several cases we can only grasp parts of the phenomenon we want to model, so we do not necessarily find a uniform mapping. Overall, what we do need is that a series of morphisms creating and maintaining the connection between the target domain and the different partial maps in a coherent way. 

It is analogous to the case recently discussed by \cite{pisano2023instrumentalist}, who argues that in the study of adaptive systems within the Free-Energy Principle framework models can ever be isomorphic to their targets. In that framework the so-called Markov-blankets play a somewhat similar role to partial maps in the modeling practice.

\subsection{The logarithmic function}

The logarithmic function is part of the standard curriculum, thus everyone is familiar with it. For this function, the following identity holds:
\begin{align*}
  \log_a(xy)&=\log_a(x)+\log_a(y)
\end{align*}
The shape of the equation is like Equation \ref{eq:compatible_operations} describing compatible operations.
Indeed, it is a prime example of a morphism.
The idea is that we map a computationally more demanding task to an easier one.
Nowadays, it is difficult to imagine, but multiplication was a demanding and error-prone task.
Addition was easier.
Logarithm turns multiplication into addition.
\gap
\begin{center}
\Morphism{numbers}{multiplication}{numbers}{addition}
\end{center}
\gap
The trick is to have the map precalculated in thick books, the logarithm tables. 
This is a pure example of how mathematics makes things easier by transferring the problem through a morphism.
Ironically, this is the point where most people ``exit'' studying mathematics.

\subsection{Chess Puzzle}

\begin{figure}
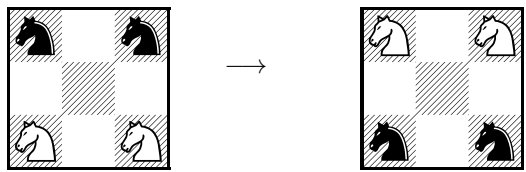

\begin{center}
 \begin{tabular}{>{\centering\arraybackslash} m{3cm} >{\centering\arraybackslash} m{1cm} >{\centering\arraybackslash} m{3cm}}
   \chessboard[maxfield=c3,addblack={Nc3,Na3},addwhite={Na1,Nc1},showmover=false,hlabel=false,vlabel=false] &
   $\longrightarrow$ &
 \chessboard[maxfield=c3,addwhite={Nc3,Na3},addblack={Na1,Nc1},showmover=false,hlabel=false,vlabel=false]
 \end{tabular}
\end{center}
\caption{Guarini's Puzzle asks whether it is possible to rearrange the knights by the usual knight's move, reaching the position on the right starting from the position on the left on a miniature 3$\times$3 chess board.  If possible, how many moves it takes to swap the black and
  white knights?}
  \label{fig:chess_puzzle}
\end{figure}

A famous chess puzzle is a good demonstration of the trick of switching representation through a morphism.
The task, see Figure \ref{fig:chess_puzzle}, is to get from a certain configuration to another by using the legal chess moves of the knight.
If the rearrangement is possible, we need to count the number of moves (\cite{levitin2011algorithmic}).

This task is difficult enough that people struggle to do it in their heads.
Even if they are given the small chess board and the pieces, they make mistake in counting the number of moves needed.
One possible cause is the unusual movement of the chess knight.
A suitable morphism should forget about that confusing detail.

\begin{figure}
\begin{center}
  \begin{tabular}{ m{3cm} m{3cm} m{4cm} }
  \begin{tabular}{|c|c|c|}
    \hline
    1 $B_1$ & 2  & 3 $B_2$\\
    \hline
    4 & 5 & 6 \\
    \hline
    7 $W_1$ & 8 & 9 $W_2$ \\
    \hline
  \end{tabular} &
  \begin{tikzpicture}
    \graph[clockwise, radius=1.2cm, n=9,nodes={draw, circle,inner sep=1pt}]
    {
        1,2,3,4,5,6,7,8,9;
        1--6--7--2--9--4--3--8--1;
    };
  
\end{tikzpicture}
  &
    \begin{tikzpicture}
    \graph[clockwise, radius=1.2cm, n=8,nodes={draw, circle,inner sep=1pt}]
    {
        1,6,7,2,9,4,3,8;
        1--6--7--2--9--4--3--8--1;
    };
    \node [above=1pt of 1]  (B1) {$B_1$};
    \node [left=1pt of 3] {$B_2$};
    \node [right=1pt of 7] {$W_1$};
    \node [below=1pt of 9] (W2) {$W_2$};
    \draw [->, very thick] (1) to[out=-30,in=30] (9);
\end{tikzpicture}
\end{tabular}
\caption{The table on the left shows how to assign numbers to the squares of the chess board and showing the initial positions ($B_1, B_2$ black, and $W_1,W_2$ white knights). The messy graph in the middle is an attempt to draw the neighbouring relationships by the knight's move, but the unnecessary morphism from the order of natural numbers ruins the usefulness of the graph. On the right we draw the graph by closely following the knight's move. The thick arrows shows the 4 moves needed for the black knight $B_1$.}
\label{fig:solution}
\end{center}

\end{figure}
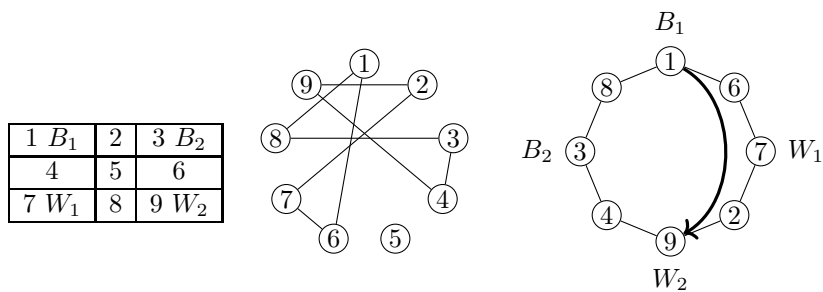

The narrative of the solution starts with abstracting away the details of the board (Fig.~\ref{fig:solution}).
We can just attach some symbol to each square.
Positive numbers might be an obvious choice.
Then, we can reduce a knight move to a relationship between two numbers.
In a teaching scenario, we can draw a graph, with numbers in a circle and the edges between them according to the jumps.
The drawing is a mess, since we assumed the numbers are nicely in order instead of the actual neighbouring relationship.
This is an example of a wrongly applied morphism.
The numbers here are just symbols, their successor relations do not transfer to the graph.
However, the graph drawn properly shows that each knight should make 4 jumps each, adding up to 16 in total.
It is just a rotation of the diagram.
Thus, by using the morphism
\gap
\begin{center}
  \Morphism{board positions}{knight moves}{nodes}{edges}
\end{center}
\gap
we turned a difficult problem into an easy one: instead of combining chess knight jumps, we simply slide them along the circle of nodes.

\section{Existing Implicit and Explicit uses of Morphisms}
\label{background}

In many fields, people already talk about morphic relations, sometimes without realizing it.
Others use only a partial definition of morphisms or do not consider them relevant to scientific or philosophical problems.
Here, we mention some salient examples.
Our thesis comes from a generalization of these cases.

\subsection{Category Theory}
\begin{quote}
 ``\ldots we are probably at least 100 years away from
a world in which one can adequately realize that category theory is everything
philosophy ever strove to be \ldots'' \cite{2022_sheaf_Rosiak}
\end{quote}
Informally, category theory can be described as an abstract theory of functions
(\cite{awodey2010category}).
As a relatively concrete example of a \emph{category}, we mention algebraic \emph{objects} of the same
kind (e.g.,
groups, rings), and the homomorphisms between them represented by \emph{arrows}.
Visually, a category is a directed graph, where associative composition is defined.
For a path, a composable sequence of arrows, there is an arrow for the corresponding composition.
In other words, whenever there is a way (no matter how long) from one object to another, there has to be a direct link.
It is rather surprising this minimal structure can serve as the model of numerous mathematical theories (\cite{spivak2017}).

Agreeing with the opening quote, we think that as category theory runs its
course, serving as an organizational meta-theory for branches of mathematics and
getting more applied elsewhere, everything we say here about
morphisms will become common knowledge.
This prediction can be justified by the progress of applied category theory (e.g., \cite{fong2019invitation,spivak2014category}).
However, there are features slowing down this process. Currently, category
theory for our purposes is
\begin{enumerate}
\item too specific to mathematics;
\item too general about morphisms.
\end{enumerate}
There is no contradiction here.
Learning category theory requires wide and mature mathematical knowledge, 
although recent books (\cite{cheng2022joy,2022_sheaf_Rosiak}) ease this situation.
Also, the exact definition of a category is geared towards algebraic objects.
The need for the identity morphism for each object is less relevant for a
general philosophical theory of explanation and representation.
On the other hand, an abstract  morphism in a category does not need to be a structure preserving map.
For instance, a partially ordered set is a category.
The arrows are the $\leq$ relations.
These arrows merely indicate the presence of the relation without being morphisms in the literal sense.

The most relevant feature of category theory here is that it provides tools from transforming one type of mathematics into another type.
It studies knowledge transferring.
Moreover, category theory has a self-similar nature: a functor is a structure preserving
map between categories, and natural transformation are morphisms between
functors.
Thus, we have a category of categories (assuming provisions against the
paradoxes of naive set theory).
Then, we can proceed with structure preserving maps between functors.
This plurality of cross-level structure preserving maps is what we envision here
for a wider context of the philosophical discourse about the interrelated notions of scientific explanation, understanding, reasoning, representation, and intelligence in general.

\subsection{Algebraic Automata Theory}

Algebraic automata theory studies computation with the tools of abstract algebra.
We represent finite state automata as \emph{semigroups} and their relationship as morphisms.  
One fundamental result, the prime decomposition theorem of finite state automata,  states that any automaton can be built from simpler components using series-parallel connections \cite{primedecomp65}.
The control information flow is unidirectional between the components, thus the decompositions are \emph{hierarchical}.
An \emph{emulation}, a surjective morphism from the original automaton into a nicer decomposition, allow us to understand original automaton better: we can see its subsystems and give coarse-grained descriptions of its dynamics.
Decompositions are not unique, thus we can understand something in several ways.

A finite state automaton (with no specified initial and accepting states for language recognition) is a \emph{discrete dynamical system}.
The above decompositions introduce \emph{coordinates} into these systems describing their dynamics on several levels.
For instance, \emph{conserved quantities} in physical systems (e.g., energy) give the top level coordinates in a hierarchical decomposition \cite{wildbook}, thus the hierarchical decomposition theory can be viewed as one basic mechanism for scientific understanding.
The theory can be extended to continuous systems \cite{KRENER1975}, and a recent reformulation on categorical level \cite{egrinagy2025sgpoid} allows further generalizations. 

\subsection{Cognitive Metaphors}

Cognitive metaphors \cite{lakoff2008metaphors} are thought processes where we
``experience and understand one kind of things in terms of another''.
Metaphor is a map between a source domain (what we know) and the target domain
(what is new for us) with the purpose of conceptualizing the target.
This is exactly the setup for a morphism, but the correspondence between the domains is not detailed. 
The authors indeed started with the idea of a mathematical mapping (described in the
2003 Afterword), but found it unsuitable for creating target domains.
Mathematical functions assume the existence of the domains, seemingly ruling out the possibility of creative metaphors.
This shortcoming is fixed by the generating sets for domains, we establish the morphisms for the generator elements only.
It is true that metaphors have limited scope, and indeed in everyday language understanding is always partial.
Chance is always there that we just ask one more question and it becomes clear that our conversion partner has different ideas.
Morphisms defined for generators guarantee morphisms for the whole domain only for mathematical objects.
Nevertheless, the mechanism is the same for partial agreements. 

A complete analysis of the metaphorical structure of Mathematics followed in \cite{lakoff2000mathematics}, but metaphors were used only as cognitive linguistic tools external to the subject, not as  mathematical objects.
This was deliberate since the main purpose was to emphasize that the philosophy of Mathematics is better start with cognitive science, rather than Mathematics itself.
Still, making the connection with category theory could have been beneficial.

Another similar missed opportunity in \cite{lakoff2008metaphors} happens when the authors distinguish between \emph{inherent} properties and \emph{interactional}
properties of entities.
They argue that our concepts are primarily defined by how they interact with each other.
As a perfect match, in category theory the Yoneda Lemma (see e.g.,\cite{awodey2010category}) states that an
object can be fully characterized by its relations to other objects, i.e., by the morphisms into and out of the object.
Despite the omission of this connection, the cognitive metaphors' overall picture of concepts cross-connected with knowledge-transferring mappings is very similar to our perspective.

\subsection{Mathematical explanation and representation}

The role of morphism in (scientific) representation plays an important role in philosophical discourse, especially in the structuralist approach, since the seminal paper of \cite{van1980scientific}. We will discuss this crucial aspect in detail in section \ref{explanation}. 

The concept of morphism also appears in the philosophical discourse of explanation (see e.g. \cite{bartels2006defending, baron2017mathematics}).
Particularly, it is mentioned explicitly in the concept in which
the explanation of physical phenomena based on mathematical models
is examined with counterfactual means (\cite{baron2017mathematics, baron2020counterfactual}). 

The basic idea is that between the two structures, the mathematical
structure and the physical structure there must be some link, such that a change 
to the mathematical structure implies a corresponding change in the physical structure,
and this link "must be held fixed as part of the counterfactual supposition
in order for the mathematical twiddle to properly ramify" (\cite{baron2017mathematics}). 

This link is nothing else, than a morphism, but the authors do not use the concept
precisely enough for it to be properly exploited. As later critics (\cite{ knowles2021platonic, kasirzadeh2023counter}) also noted,
one of the main problems with the approach is that it tries to impose changes
in the mathematical structure on the physical structure. 

At the same time, the critics also make a very typical mistake in that they necessarily
mean homomorphism or even isomorphism by morphism (\cite{pero2016varieties}), which
``would not survive certain changes to the mathematical domain'' (\cite{knowles2021platonic}). 
We will come back to this issue and clarify the difference soon in chapter \ref{definition}.

\subsection{The Thousand Brains Theory of Intelligence}
\label{sect:1000brains}
\begin{quote}
  ``Grid cells in the neocortex suggests that all knowledge is
  learned and stored in the context of locations and location spaces
  and that “thinking” is movement through those location spaces.'' \cite{2019hawkinsarticle}
 \end{quote}
We described the process of thinking as a composable sequence of events and understanding is matching the paths through the space of ideas.
This might be literally true according to the Thousand Brains Theory \cite{hawkins2021thousand}.
It defines intelligence is the ability of creating predictive models of the world.
The novelty of the theory lies in positing that several models of different sensory inputs exist for a single object, hence the name thousand brains, instead of having one model for one thing.
Despite the multitude of models, we have stable, singular perceptions.
Connections between the models and a voting mechanism ensure unambiguity.
This mechanism is like a low-level version of understanding as finding structural matches.
Conversely, our theory of understanding brings this voting mechanism for perception to a higher level of thinking.

The representations are spatial for all models, using reference frames.
These are like coordinate systems but not as simple as Euclidean space.
According to the hypothesis, the grid cells in the neocortex used for representing locations are reused for abstract thinking \cite{2019hawkinsarticle}.
In this view, thinking is movement.

\section{Explanation, representation, understanding and communication}
\label{explanation}

How can we describe explanation and understanding using morphisms?
By dictionary definition, an explanation is a set of statements used to answer questions. Typically, the questions have the form ``Why \ldots?''.
This traditional definition is a particular case: a sequence of statements connected by logical inference leads to the \emph{explanandum}.
We can abstract away the logical inference, as proceeding from one thought to the other can be just recalling memories, or finding associations.
At this stage, we consider only what mechanism counts as an explanation, but not the conditions that make an explanation correct.
What remains after this abstraction is a path in a directed graph leading to some representation of the event to be explained.
If we have a morphism to find a compatible path in the problem domain, then we have an explanation.
Therefore, \emph{explanation is pathfinding in the easier domain of a morphic relation.}

The underlying assumption is that thinking is a sequential process in the context of explanation.
This is not to deny other thinking processes.
The brain is massively parallel, and sudden creative insights are possible even in mathematics.
However, for the modes of thinking (e.g., scientific, mathematical), where we have the possibility of communication, the sequential process is a requirement. 

Understanding is a special case: \emph{an explanation given to myself}.
It is about bringing my thought sequence into a morphic relation with the explanandum.
The condition of compatible operations gives choice.
I can follow steps in my thinking, only mapping the last one to the phenomenon.
Or, I can translate each step in my thinking and carry out the actions or make the observations.
At the end, I should get the same result.

Planning follows the same pattern.
We can solve the problem of reaching the bananas by swinging the sticks or using other tools to see what works.
Alternatively, we can also imagine sequences of actions and evaluate them, then carry out the best action only.

\subsection{Use cases of morphisms}
The underlying idea of understanding is simple, pathfinding in the easier domain and using a morphism. 
However, there are several use cases of this simple motif.

\begin{itemize}
\item Using an existing morphism for \emph{predicting} what would happen in the target
  domain based on what we can do in the source domain (\emph{exploitation}, using the
  map to find a café);
\item realizing the presence of an existing morphism (\emph{discovery}, figuring out
  where we are on the map by associating immediate landmarks to map features);
\item constructing a morphism (\emph{creative} mode, the chess puzzle example);
\item \emph{communication}, teaching.
\end{itemize}



When two people, A and B want to come to an agreement on a problem, they must achieve a correspondence between their existing thought structures. 
When A wants to explain something to B, meanwhile B intends to understand the explanation of A, then A must help B to develop a thought structure similar to the one A already possessed.
Interoperability between these two structures of A and B can be provided by various ways, but in order to have a chance that they \emph{really} mean the same thing under a certain problem or explanation, this mapping must preserve the effect of the (logical or other type of) operations taking place in the thought structures.
Such a mapping is nothing else than a morphism.
That is, when two people agreed on something, or one successfully explained something to the other,
what they finally have is a morphic relationship between their thought structures.

From this point of view, agreement or explanation is based on the same foundations. 
The difference is that in the case of agreement, the key moment is finding the morphism between the two solid, developed thought structures of the two persons, A and B, while in the case of explanation, the explaining person (A) wants to create a thought structure in B equivalent (up to a morphism) to his or her own thought structure. 
For the latter, A can apply elements of a realized basic morphism already established by the basic thought elements existing in both persons (for example coming from the structure of the common language, or from the similar education experience), but A must build the thought structure in B parallel to the morphism related to the given area, constantly taking care to preserve the morphic relationship between the two (the existing and the currently being developed) structures.

One test for shared understanding could be that two persons tell the same explanation swapping the speaker role frequently (like in team presentations), since in a morphism does not matter where we make the next move.
Of course, even this is not an absolute guarantee of having the same idea, as different thought structures can map partially.
However, without it, agreement or explanation will almost certainly fail.

Finally, it is worth noting that some authors believe that our first item about using an existing morphism between the model and the target domain of the reality for \emph{predicting} what would happen in that domain can further be sharpened to using this morphism for \emph{affecting} what would happen in the target domain (\cite{tee2019constructing}). Although this is beyond our scope of this paper, it is clear that this affection could happen only through a general morphic relationship between the two domains.

\subsection{Maps revisited}

Now let us extend the idea of mapping and consider the map as a representation or explanation (depending on its usage) of a domain. The cartographic map, as we have seen in the previous section, represents a slice of the world to us based on a morphism: if you go this way, you will get there; if that's what you're looking for, you'll find it here. And this will also happen in reality, with smaller but still tolerable differences due to the restrictions mentioned above.

If the map is an explanation, then different maps, different morphisms provide different explanations for the same phenomenon. Which one do we prefer? One possibility would be to choose a morphism that is an explanation which perfectly explains the world and its structure we want to understand from one specific aspect, while strongly distorting it from another aspect. In most cases, however, we prefer an explanation that ``roughly'' explains the phenomenon we want to learn about well (well, but not perfectly) from all aspects. This is a ``generally distorted'' explanation, a morphism through which we get an overall acceptable picture of the phenomenon.

When we want to explain a thought structure or gain knowledge about it by understanding an explanation of somebody else, it can happen that it is theoretically or practically impossible to map this structure as a whole, but some parts of it can still be explained well, which is analogous to the situation when there are only local maps of a world. Another explanation can explain some other aspects of this thought structure, but these two different explanations should be consistent in the same sense what we have seen in the case of local mappings.

\subsection{Games}

When thinking about games the meaningless nature of the rules is often mentioned.
The formalist view of mathematics, starting with Thomae and Frege (\cite{frege1903}) especially likes to make this claim, using the
chess analogy to show that arithmetic symbol manipulation has no further meaning (for a recent overview see \cite{lawrence2023frege}).
We argue that this is incorrect, and if the rules were meaningless, then no one
would play the game.
The meaning of games is multifaceted and defined through several morphisms.

One of the motivations for playing traditional board games is the desire to win.
Thus, we have the morphism \morphism{Board Positions}{Emotional States}  for most
human players.
The operations for the domain of board positions are the game moves (like in the game
tree), and  the compatible operations are movements in the multidimensional
space of emotions.
The target domain is of course not defined quantitatively, but the compatibility is clear.
When observing a player, one can concentrate on the board, and by seeing advantageous moves, we infer a happier state of mind.
Or, if we look at the facial expressions of the player, we can predict how well the game is going for her.

Another human goal is to understand the game, i.e., giving explanations for winning and losing in the game.
A master player game analysis provides a \morphism{Game
  Record}{Narrative} morphism.
The narrative has a causal structure.
If this `storyline' is compatible with the moves, then we analyzed the game successfully.
From a player perspective, we can talk about plans instead of narratives.

In machine learning the morphism \morphism{Board Positions}{Winning Probabilities} is the fundamental tool.
The target domain is one-dimensional, the $[0,1]$ interval of real numbers.
This reduces information heavily, and is often criticized for not providing understanding.

All these morphisms are relative to strength of the player (both human and AI).
The same move can map to totally different emotional states (even after sign
correction).
Different AI engines can give different winning probabilities depending on their
network size and training time and data set.

For a pure skill two-player zero-sum game there is ground truth defined by the
fully evaluated game tree.
Then we have the \morphism{Board Position}{Guaranteed Score} morphism for solved
games.
For complex enough games (like Chess and Go) this is computationally unreachable.
A non-perfect player tries to approximate this ground truth morphism.

\subsection{ Scientific Understanding}

As we have mentioned, the role of morphism plays crucial role in the structuralist approach of (scientific) representation and explanation (\cite{van1980scientific}). 
While the basic idea is accepted by most researchers, criticism appears in terms of its type and usability. The basis of the negative opinions is that the articles promoting this idea often interpret morphism too narrowly, which in this sense leads to legitimate criticism. Van Fraasen himself could only imagine isomorphism as scientific representation (\cite{van1993vicious}, \cite{van1994interpretation}), but later researchers who wanted to refine this principle faced the same issue as van Fraasen. 
Partial isomorphism  (\cite{bueno1997empirical}; \cite{french2003model}); homomorphism (\cite{lloyd2021structure}; \cite{bartels2006defending}), or $\Delta/\Psi$-morphism (\cite{swoyer1991structural}) have been proposed to weaken strict conditions, but these efforts always result in some narrowing of the general morphism, while we claim that the essence of representation and explanation is precisely the existence of the (dynamic) morphism, not the fulfillment of various extra properties. 

General morphism allows us the necessary flexibility, and with any narrowing, all we achieve is to find counterexamples where they do not work. While the main thing is to ensure the existence of the structure and its preservation to a certain degree, no matter how weak these conditions are. 
Moreover, \cite{boesch2021scientific} puts forward the very valid argument that too strong criteria and too accurate (e.g. isomorphic) relationships can make representation impossible or empty.
It is no coincidence that structures that are isomorphic to each other are often regarded as identical in pure mathematics, either treated as ``abuse of notation'', or as a new foundational axiom (univalence axiom \cite{awodey2013univalence}).
 A perfect correspondence has no relevant explanatory power and scientific challenge and inspiration beyond the fact that the two structures are the same. As \cite{boesch2021scientific} argues, it is dissimilarities, or weakened morphisms that inspire scientific discovery.

The existence of the morphism is therefore sufficient. On the other hand, we consider the existence of the morphism to be necessary in all such situations. In the case of (scientific) model creation and representation, several researchers criticize the various special morphisms, which are excellently summarized by \cite{suarez2003scientific} and \cite{boesch2021scientific}. In these criticisms, it is also suggested that in fact the representation does not necessarily require any morphism, because there is not necessarily an operation or structure that we want to preserve.
We may think about a representation without a structure through a simple example. Indeed, integers, or the stations on the London Underground map are interesting in themselves. However, when we look at them as a model, a representation, and want to explain something with their help, we cannot avoid using some relation or operation that makes them a structure. How do I get to this part of town, or how can I set up a secret code with the help of primes - these are actions where the mere set is enriched with relations and operations and becomes a structure. We must preserve this structure when using the model and representation in a meaningful way.

Similar arguments support the necessity of having a morphism related to biomedical research by \cite{varga2024scientific}. 
As it is stated, ``understanding does not only require the possession of a theory or model, but also the skill or ability to use it to discern the causal relationship involved.'' 
Without a morphism, the theory is somewhat ``disabled'', since it cannot help finding a causal relationship, the essence of understanding. The morphism is therefore required by nothing else than the meaningful use of the representation.

\subsection{(The lack of) morphisms of mindmaps}

In recent years, the use of software specifically designed to visualize thoughts and explanations has become widespread. This procedure is known by various names, such as concept mapping, mind mapping and argument mapping (\cite{davies2011concept}). 
This approach has become particularly popular in education and scientific research.

The basic concept is that users can represent or manipulate relationships among thought units in a diagram, practically in a directed graph with weighted connections. 
These graphs can be very suitable for explaining a concept, but in the literature and during use, people frequently ignore a crucial aspect, and this is morphism.

The mapping of a more complex thought structure by different people is likely to yield different graphs. 
It may happen that these graphs are not even isomorphic, or if they are, then the thought elements appearing in individual vertices or the individual connections do not represent the same thought items in two different mind mappings. 
Separately, both can be a valid representation of the thoughts of the creator, but without structure preserving correspondence between the two maps the mutual understanding is practically hopeless.

\subsection{The issue of (un)explainable AI}

In recent years, explainability has been one of the most important issues of Artificial Intelligence, both ethically and epistemologically. 
AI techniques, especially with the appearance of deep learning and generative AI, are becoming increasingly opaque, and even if we understand the underlying mathematical principles of such models they still lack explicit declarative knowledge (\cite{holzinger2018machine}). 
The basic problem is the incredible increase in the amount of data and the complexity of the structure of our AI models: we should find a morphism from these structures and decisions appearing in many dimensional representation spaces into dimensions that can be interpreted by human intelligence - however, this is practically hopeless with such a reduction in the number of dimensions.

This does not contradict the fact that we can handle and interpret abstract spaces of arbitrary dimensions in a mathematical sense. 
In the latter case, the structure is historically created as a generalization of a low-dimensional structure.
For example, $\R^3$ has been existed and interpreted by humans for centuries, and the thought structure of $\R^n$ for arbitrary large $n$ was built on that model. 
However, in the case of AI, ``we are challenged with data of arbitrarily high dimensions, and within such data, relevant structural patterns and/or temporal patterns (``knowledge'') are often hidden'' (\cite{holzinger2018machine}), and cannot be reduced or mapped to lower dimensional spaces – and even less likely to be explained by any morphism.

The latest contributions, such as the thought-provoking approach of \cite{buchholz2023means}, also implicitly look for some morphism, e.g. between the internal state space of the AI and some humanly interpretable domain, even by interposing some intermediate domain, e.g., by heatmapping.

\section{Conclusion}
\label{conclusion}

In this paper, we argued that morphism is a fundamental element of all our intelligent activities, and a universal mechanism for any understanding, representation, reasoning and explanation, for human as well as for artificial intelligence.
Morphism is necessary and sufficient to establish a meaningful connection between two intelligent agents. 
Morphism is necessary and sufficient to create a meaningful representation of a target domain in another domain. We have provided several examples of how this morphism appears in different fields. 
It has also been emphasized that any narrowing of the concept of morphism inevitably leads to obstacles and counterexamples.

We can also consider this argument from a wider perspective.
We live in a regular enough world to have executable morphic models useful for survival and, eventually, for understanding it.
Being executable means the models are dynamic; we can move through some space of the ideas.
Being morphic means that operations in the different domains are compatible; we can choose the domain to work in.
We can think about composing actions or perform them in the world.
For morphic relations, the result of thinking stays in correspondence with the result of the action sequence.

These cognitive models of how our everyday world works are specific examples.
However, evolutionary and historically speaking, they might have been the first morphisms.
As a second step, morphic relations appeared between the models, enabling communication
and understanding by finding compatible operations in thinking.
To be successful in a hunter-gatherer tribe, members needed to have reliable morphic relations between their morphic models of their environment: the weather, the prey, and hunting weapons.

We can define Mathematics as the limit of this process: using and creating morphisms in overdrive.
The mathematical examples of morphisms (e.g., group and graph homomorphisms, diffeomorphisms and topological maps, functors between categories) are distilled manifestations of \emph{the mechanism of compatible operations between various domains - a necessary ingredient for any intelligence}.
However, all other fields of knowledge can benefit from the explicit description of compatible operations and their domains.

\section*{Acknowledgements}
This project was partially funded by the Kakenhi grant
22K00015 by the Japan Society for the Promotion of Science (JSPS), titled `On progressing human understanding in the shadow of superhuman
deep learning artificial intelligence entities' (Grant-in-Aid for Scientific
Research type C, \url{https://kaken.nii.ac.jp/grant/KAKENHI-PROJECT-22K00015/}).

\printbibliography

\end{document}